\documentclass{amsart}

\makeatletter
\@addtoreset{equation}{section}\makeatother
\newtheorem{theo}{Theorem}[section]
\newtheorem{lem}[theo]{Lemma}

\theoremstyle{remark} \newtheorem{remark}[theo]{Remark}
\newtheorem{example}[theo]{Example}

\usepackage{amssymb}
\usepackage{amscd}
\usepackage{verbatim,ifthen,epsfig}

\usepackage{amsmath,graphics}
\usepackage{latexsym}
\usepackage{eucal}

\usepackage{latexsym}

\def\RR{\mathbb{R}}

\def\ang#1{{\langle #1 \rangle }}

\def\Id{\operatorname{Id}}

\def\ang#1{{\langle #1 \rangle }}

\newcommand\spec{{\operatorname{spec}}}
\newcommand\umin{u_{\operatorname{min}}}
\newcommand\vmin{v_{\operatorname{min}}}
\renewcommand\tt{\tilde t_{\Id}}
\newcommand\ttt{\tilde t_F}
\newcommand\tttm{\tilde t_{F_\mu}}
\newcommand\Deltab{\Delta_{\partial \Omega}}
\newcommand\Et{\tilde E}
\newcommand\Omegabar{\overline{\Omega}}
\newcommand\Omegab{\partial \Omega}

\newcommand{\be}{\begin{equation}}
\newcommand{\ee}{\end{equation}}

\begin{document}
\title[Neumann eigenfunctions at the boundary] {Estimates on Neumann eigenfunctions at the boundary, and the ``Method of Particular Solutions'' for computing them}

\author{A. H.  Barnett}
\address{Department of Mathematics, Dartmouth College, Hanover, NH, 03755, USA}
\email{ahb@math.dartmouth.edu}
\author{Andrew Hassell} %
\address{Department of Mathematics, Australian National University, Canberra ACT 0200, AUSTRALIA}
\email{Andrew.Hassell@anu.edu.au}

\begin{abstract}
We consider the {\em method of particular solutions} for numerically computing eigenvalues and eigenfunctions of the Laplacian on a smooth, bound\-ed domain $\Omega$ in $\RR^n$ with either Dirichlet or Neumann boundary conditions.
This method constructs approximate eigenvalues $E$, and approximate
eigenfunctions $u$ that satisfy $\Delta u=Eu$ in $\Omega$, but not the
exact boundary condition.
An inclusion bound is then an estimate on the distance of $E$ from the actual spectrum of the Laplacian, in terms of
(boundary data of) $u$.
We prove operator norm estimates on certain operators on $L^2(\partial \Omega)$ constructed from the boundary values of the true eigenfunctions, and show that these estimates lead to sharp inclusion bounds
in the sense that their scaling with $E$ is optimal. This
is advantageous for the accurate computation of large eigenvalues.
The Dirichlet case can be treated using elementary arguments and will appear in \cite{bnds}, while the Neumann case seems to require much more sophisticated technology.
We include preliminary numerical examples for the Neumann case.
\end{abstract}
\maketitle

\section{Introduction}
In this paper we consider Laplace eigenfunctions on a smooth bounded domain $\Omega \subset \RR^n$. As is well known, the positive Laplacian\footnote{Note that our sign convention is opposite to that of \cite{bnds}},
$$
\Delta = - \sum_{i=1}^n \frac{\partial^2}{\partial x_i^2},
$$
with domain either $H^2(\Omega) \cap H^1_0(\Omega)$ or $$
\{ u \in H^2(\Omega) \mid d_n u |_{\partial \Omega} = 0 \}
$$
is self-adjoint. Here and below, $d_n$ denotes the directional derivative with respect to the outward unit normal vector at $\partial \Omega$. These are known as the Laplacian with Dirichlet, resp. Neumann, boundary conditions and will be denoted $\Delta_D$, resp. $\Delta_N$. In either case, there is an orthonormal basis of $L^2(\Omega)$ consisting of real eigenfunctions. We denote an orthonormal basis of Dirichlet eigenfunctions $u_j$, $j = 1 \dots \infty$, and an orthonormal basis of Neumann eigenfunctions $v_j$, $j = 1 \dots \infty$, with Dirichlet/Neumann eigenvalues $E_j = \lambda_j^2$, resp. $\Et_j = \mu_j^2$. Thus $u_j$, $v_j$ satisfy a \emph{Helmholtz equation}
$$
\Delta u_j = E_j u_j = \lambda_j^2 u_j , \quad \mbox{ or } \quad
\Delta v_j = \Et_j v_j = \mu_j^2 v_j,
$$
with boundary condition
$$
u_j |_{\partial \Omega} = 0, \quad \mbox{ or } \quad
d_n v_j|_{\partial \Omega} = 0.
$$
We will denote the spectrum of $\Delta_D$, resp. $\Delta_N$ by $\spec_D$, resp. $\spec_N$. Also, we will denote the normal derivative $d_n u_j$ at $\partial \Omega$ by $\psi_j$, and the restriction of $v_j$ to $\partial \Omega$ by $w_j$. Thus $\psi_j, w_j$ are functions on $\partial \Omega$, which we will refer to as the boundary traces of eigenfunctions $u_j$, resp. $v_j$.

\begin{figure} \epsfig{file=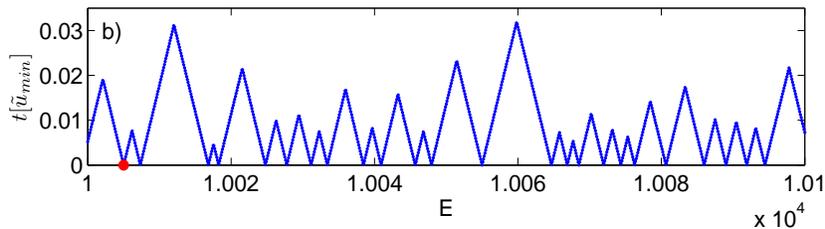,width=11cm,height=3cm}
\caption{\small Graph of the minimum value of $t[u]$ achievable at each $E$,
for the Dirichlet case; see \eqref{ib}.
The smooth domain $\Omega\subset\mathbb{R}^2$ is shown in Fig.~\ref{f:mode}.
$u$ is restricted to lie in the span of a finite number of
numerical basis functions satisfying the Helmholtz equation at each $E$,
and the minimization uses dense linear algebra \cite{bnds}.
A small value (e.g. at the dot shown)
implies closeness to a Dirichlet eigenvalue.
}
\label{f:tE}
\end{figure}

The Method of Particular Solutions \cite{mps,incl} is a numerical method for finding eigenvalues and eigenfunctions of the Laplacian on a Euclidean domain. In the case of the Dirichlet boundary condition, the method consists of choosing an energy (positive real number) $E$, and looking for the solution $u$ to the Helmholtz equation $(\Delta -E) u = 0$ that comes closest to satisfying the boundary condition, in the sense that it minimizes (or approximately minimizes) the $L^2$ norm of the boundary trace of $u$  (in $L^2(\partial \Omega)$).
In practice $u$ is restricted to a sufficiently large numerical subspace
and the minimization is done via dense linear algebra
(a generalized eigenvalue or singular value problem \cite{incl,bnds}).
We then think of this minimum $L^2$ norm on the boundary as a function of $E$
(see Fig.~\ref{f:tE}), and numerically try to find the (near) zeros of this function as we move along the $E$-axis. Clearly, if we find a Helmholtz $u$ with $\| u \|_{L^2(\partial \Omega)}$, then $E$ is a Dirichlet eigenvalue, and $u$ is a Dirichlet eigenfunction. We would expect, therefore, that if $\| u \|_{L^2(\partial \Omega)}$ is very small, then $E$ is close to a Dirichlet eigenvalue, and $u$ is close to (i.e. makes a small angle with) the corresponding eigenspace. 

An \emph{inclusion bound} is a quantitative estimate of this form, taking the form (for eigenvalues)
\begin{equation}
d(E, \spec_D) \leq C  E^\alpha t[u], \quad \mbox{ where } \quad t[u] := \frac{ \| u \|_{L^2(\partial \Omega)}}{\| u \|_{L^2(\Omega)}}
\label{ib}\end{equation}
for some constant $C$ independent of $E$ and exponent $\alpha$. This tells us how small we must make the `tension' $t[u]$ in order to achieve any desired accuracy in our numerically computed eigenvalue. We are primarily interested in high energy estimates, i.e. $E$ large, so our goal is to obtain such an estimate which is sharp as $E \to \infty$, that is, with the smallest possible $\alpha$. (Ideally, for practical applications, we also want a constant $C$ that is small and computable, that is, expressed in terms of geometric quantities such as the measure, surface measure, inradius, etc., of our domain $\Omega$. But we will not address that issue here.)

This article is a report on completed work on the Dirichlet case \cite{bnds}, and an announcement of work in progress on the Neumann case. Full details will appear elsewhere.

\section{Dirichlet boundary condition}
We begin by proving  that there are upper and lower bounds
\begin{equation}
C^{-1} \lambda_j \leq \| \psi_j \|_{L^2(\partial \Omega)} \leq C\lambda_j \label{upperlower}\end{equation}
where $C$ depends only on $\Omega$. These are classical and well-known estimates, but we give the proof since it follows from a calculation that we need later on anyway.
The proof is via a  Rellich-type identity, involving the commutator of $\Delta$ with a suitably chosen vector field $V$. 
The basic calculation is 
$$
\ang{u, [\Delta, V]u} = \int_\Omega \Big( ((\Delta - \lambda^2) u) (V u) - u (V(\Delta - \lambda^2) u)  \Big) $$
\begin{equation} + \int_{\partial \Omega} \Big(  (d_n u) (Vu)  - u (d_n (Vu)) \Big) .
\label{basic}\end{equation}
If $u = u_j$ is a Dirichlet eigenfunction with eigenvalue $\lambda_j^2$, then three of the terms on the RHS vanish, and we obtain
$$
\ang{u, [\Delta, V]u} = \int_{\partial \Omega}   (d_n u)(Vu)   .
$$

If we choose $V$ so that, at the boundary, it is equal to the exterior unit normal, then the RHS is precisely $\| \psi_j \|^2$.
(If not indicated, norms will be assumed to be $L^2$ norms.)
The left hand side is $\ang{u, Q u}$ where $Q$ is a second order differential operator and is $O(\lambda_j^2)$, yielding the upper bound $\| \psi_j \|^2 = O(\lambda_j^2)$.
On the other hand, if we take $V$ to be the vector field $\sum_i x_i \partial_{x_i}$, then $[\Delta, V] = 2\Delta$. Then the LHS is exactly equal to $2\lambda_j^2$, while the RHS is
no bigger than $(\max_{\partial \Omega} |x|) \| \psi_j \|^2$, yielding the lower bound $\lambda_j^2 = O(\| \psi_j \|^2)$. This lower bound is due to Rellich \cite{rellich}.

It turns out that there is a very useful generalization of the upper bound in \eqref{upperlower}, proved recently by the authors,
that applies to a whole $O(1)$ frequency window:

\begin{theo}\label{main} Let $\Omega \subset \RR^n$ be a smooth bounded domain and let $\psi_i$ be defined as above. Then the operator norm of the finite rank operator
\begin{equation}
\sum_{\lambda_i \in [\lambda, \lambda + 1]} \psi_i \ang{\psi_i, \cdot} \quad : L^2(\partial \Omega) \to L^2(\partial \Omega)
\label{T}\end{equation}
is bounded by $C\lambda^2$, where $C$ depends only on $\Omega$. \end{theo}

Remarks:

\begin{itemize}

\item This is quite a strong estimate, since by the lower bound of \eqref{upperlower} there is a lower bound of the form $c \lambda^2$ on the operator norm of any \emph{one} term in the sum.

\item This is closely related to the phenomenon of `quasi-orthogonality' of $\psi_i$ and $\psi_j$, when $|\lambda_i - \lambda_j|$ is small. Indeed, this estimate implies that when $|\lambda_i - \lambda_j| \leq 1$, then the inner product $\langle \psi_i, \psi_j \rangle$ is usually small compared with $\lambda^2$. See Barnett \cite{que}.

\item This is also closely related to an identity of B\"acker, F\"urstberger, Schubert and Steiner \cite{backerbdry}. \end{itemize}

Theorem~\ref{main} is proved as follows: first, we prove the upper bound $\| d_n u \|_{L^2(\partial \Omega)} \leq C \lambda \| u \|_{L^2(\Omega)}$ is valid not just for eigenfunctions, but for approximate  eigenfunctions $u \in \operatorname{dom} \Delta_D$ such that
$$
\| (\Delta - \lambda^2) u \|_{L^2(\Omega)}=O(\lambda).
$$
In fact, the proof is almost unchanged: we use \eqref{basic} again. Now the term $(\Delta - \lambda^2) u (Vu)$ is no longer zero, but by assumption $(\Delta - \lambda^2) u $ is $O(\lambda)$ in $L^2(\Omega)$, and also $Vu$ is $O(\lambda)$ in $L^2(\Omega)$, so by Cauchy-Schwarz this term is $O(\lambda^2)$. We treat the term $u (V(\Delta - \lambda^2) u)$ similarly, after first integrating the vector field $V$ by parts (which produces no boundary term since $u$ vanishes at $\partial \Omega$). The rest of the argument runs as above. This was also noticed by Xu \cite{Xu}.
 Notice that this condition applies in particular to a spectral cluster, that is, for $u \in \operatorname{range} E_{[\lambda, \lambda + 1]}(\sqrt{\Delta_D})$. We then use a $TT^*$ argument:
We define an operator $T$ from $\operatorname{range} E_{[\lambda, \lambda + 1]}(\sqrt{\Delta_D})$ to $L^2(\partial \Omega)$ by $Tu = d_n u |_{\partial \Omega}$. We can express $T$ in terms of the eigenfunctions by $$
Tu = \sum_{\lambda_i \in [\lambda, \lambda + 1]} \ang{u, u_i} \psi_i.
$$
That is, it is just the normal derivative of the element $\sum a_i u_i$ of this spectral subspace. Then, as we have just shown, 
$$
\| T \| \leq C\lambda.
$$
It follows that $T T^* : L^2(\partial \Omega) \to L^2(\partial \Omega)$ has operator norm bounded by $C^2 \lambda^2$. But $T T^*$ is precisely the operator \eqref{T} appearing in the statement of the theorem.

\section{Dirichlet inclusion bound}\label{sec:Dinc}

As mentioned above, a Dirichlet inclusion bound is an estimate on the distance from $E$ to $\spec_D$ in terms of $t[u]$ defined in \eqref{ib} of a Helmholtz solution $u$. A  classical result along these lines is the  Moler-Payne inclusion bound \cite{molerpayne}. This says that
$$
d(E, \spec_D) \leq C E t[u], \quad  (\Delta - E) u = 0,
$$
where $t[u]$ is defined in \eqref{ib} and $C$ depends only on $\Omega$. The proof in \cite{molerpayne} uses very little about the Dirichlet problem in particular.

Recently, Barnett \cite{incl}, followed by the authors
\cite{bnds}, improved this bound by a factor of $\sqrt{E}$:

\begin{theo}\label{MPS-Dir} There exist constants $c, C$ depending only on $\Omega$ such that the following holds.  Let $u$ be any nonzero solution of $(\Delta - E) u = 0$ in $C^\infty(\Omegabar)$, and let $\umin$ be the Helmholtz solution minimizing $t[u]$. Then
$$
c \sqrt{E}  t[\umin] \leq d(E, \spec_D) \leq C \sqrt{E} t[u].
$$
\end{theo}

\begin{remark} There is always a Helmholtz solution that minimizes $t[u]$; see \cite{bnds}.
\end{remark}

\begin{proof} The result is trivial if $E \in \spec_D \Delta$. Suppose that $E$ is not an eigenvalue, and consider the map $Z(E)$ that takes $f \in L^2(\partial \Omega)$ to the (unique) solution $u$ of the equation
$$
(\Delta - E) u = 0, u |_{\partial \Omega} = f.
$$
The $u$ that minimizes $t[u]$ then maximizes $\| u \|_{L^2(\Omega)}$ given $\| u \|_{L^2(\partial \Omega)}$. So \be
(\min_u t[u])^{-1} = \| Z(E) \| \implies (\min_u t[u])^{-2} = \| A(E) \| ,
\label{e:tinvsq}
\ee
where the operator $A(E) : L^2(\partial  \Omega) \to L^2(\partial  \Omega)$ is defined by $A(E) =   Z(E)^* Z(E)$. We claim that $A(E)$  has the expression \cite{incl}
\begin{equation}
A(E)  = \sum_j \frac{\psi_j \ang{\psi_j, \cdot}}{(E - E_j)^2}.
\label{Aexpression}\end{equation}
Here and in the remainder of the article, we take the symbol $\displaystyle \min_u$
to mean the minimum over the space of Helmholtz solutions.

To prove \eqref{Aexpression}, we show that $Z(E)$ has the expression
\begin{equation}
Z(E)f = \sum_i \frac{\ang{f, \psi_i} u_i}{E - E_i}
\label{Zexp}\end{equation}
from which \eqref{Aexpression} follows immediately.
To express $Z(\lambda)$, suppose $f$ is given and $u = Z(E) f$. We write $u = \sum a_i u_i$ as a linear combination of Dirichlet eigenfunctions. Then,
using the Helmholtz formula and Green's identities,
$$\begin{gathered}
a_i = \ang{u, u_i} = \frac1{E - E_i} \int_\Omega \big( (\Delta u) u_i - u(\Delta u_i) \big)
\\
= \frac1{E - E_i} \int_{\partial \Omega} \big( u (d_n u_i) - (d_n u) u_i  \big)
= \frac1{E - E_i} \int_{\partial \Omega} f \psi_i
\end{gathered}$$
which proves \eqref{Zexp}.

The lower bound in Theorem~\ref{MPS-Dir} is easy to prove: we note that $A(E)$ is a sum of positive operators in \eqref{Aexpression}, so the operator norm of $A(E)$ is bounded below by the operator norm of any one summand. So, using the upper bound in \eqref{upperlower},
$$
\big\| A(E) \big\| \;\geq\; \Bigl\| \frac{ \psi_j \ang{\psi_j, \cdot}}{(E - E_j)^2} \Bigr\| \;\geq\; \frac{cE}{d(E, \spec_D)^2},
$$
where $E_j$ is the closest eigenfrequency to $E$. Since $(\min t[u])^{-2} = \| A(E) \|$, this proves the lower bound.

To prove the upper bound, we use Theorem~\ref{main}. We need to show that \begin{equation}
\| A(E) \| \leq \frac{CE}{d(E, \spec_D)^2}.
\label{want}\end{equation}

To do this, we break up the sum \eqref{Aexpression} 
into the `close' eigenfrequencies in the interval $[\lambda - 1, \lambda + 1]$ and the rest. The estimate \eqref{want} for the close eigenfrequencies is immediate from Theorem~\ref{main}.

For the far eigenfrequencies, we first treat those that lie in the interval $[\lambda/2, 2\lambda]$. These can be broken up into frequency windows of width $1$, which are distance $1, 2, 3$, etc from the chosen frequency $\lambda$. Then, in \eqref{Aexpression}, the numerator for each window have operator norm bounded by $E$, and the denominator is $n^2 E$. Since $\sum n^{-2}$ is finite, the contribution from these eigenvalues is $O(1)$.

For the the eigenfrequencies not lying in $[\lambda/2, 2\lambda]$, it is not hard to see that the contribution is only $O(E^{-1/2})$, since we get a factor $E^2$ in the denominator. So the far eigenfrequencies altogether only contribute only $O(1)$ to the operator norm of $A(E)$.

Finally we observe that $$
\frac{CE}{d(E, \spec_D)^2} \geq C
$$
for every $E$, since the distance from $E$ to the spectrum can be at most $\sim \sqrt{E}$ (this follows by considering an approximate eigenfunction supported in a ball contained in $\Omega$). Thus the contribution from the near eigenvalues dominates, and we see that
\eqref{want} is true for all $E$, completing the proof. \end{proof}

Similar reasoning gives a bound for the distance between $u$ and the closest eigenfunction:

\begin{theo} There is a constant
$C$ depending only on $\Omega$, such that the following holds.
Let $E>1$, let $E_j$ be the eigenvalue nearest to $E$,
  and let $E_k$ the next nearest distinct eigenvalue.
Suppose $u$ is a solution of
$(\Delta+E)u=0$ in $C^\infty(\Omega)$ with $\|u\|_{L^2(\Omega)}=1$,
and let $\hat{u}_j$ be the projection of $u$ onto the $E_j$ eigenspace.
Then,
\begin{equation}
\|u - \hat{u}_j\|_{L^2(\Omega)} \; \le \; C\frac{\sqrt{E}\,t[u]}{|E-E_k|}
~.
\label{e:e}
\end{equation}
\label{t:e}
\end{theo}

For the proof, see \cite{bnds}.

\section{Neumann boundary condition}

We next consider the method of particular solutions for computing Neumann eigenvalues and eigenfunctions. The Neumann boundary condition is $d_n v |_{\partial \Omega} = 0$. It seems natural to minimize  (cf. \eqref{T})
$$
\tt[v] = \frac{ \| d_n v \|_{L^2(\partial \Omega)} }{ \| v \|_{L^2(\Omega)} },
$$
over nontrivial solutions $v$ of $(\Delta - \Et) v = 0$, since $\tt[v] = 0$ implies that $\Et$ is a Neumann eigenvalue and $v$ a Neumann eigenfunction. Notice, though, that we could equally well minimize the quantity
$$
\ttt[v] = \frac{ \| F(d_n v) \|_{L^2(\partial \Omega)} }{ \| v \|_{L^2(\Omega)} },
$$
for any invertible operator $F$ on $L^2(\partial \Omega)$. It turns out that there is an essentially optimal choice of $F$ (which depends on $\Et$), which is \textbf{not} the identity. Indeed the main point of this article is to determine this optimal $F$.

The form of $F$ is suggested by the local Weyl law for boundary values of eigenfunctions. This law \cite{gerard93,hassell} says that the boundary traces of eigenfunctions are, on the average,  distributed in phase space $T^*(\partial \Omega)$ according to \begin{equation}\begin{gathered}
c (1 - |\eta|^2)^{1/4} 1_{\{ |\eta| \leq 1 \} } \text{ (Dirichlet),} \\
\phantom{aa} \tilde c (1 - |\eta|^2)^{-1/4} 1_{\{ |\eta| \leq 1 \} } \text{ (Neumann)}
\end{gathered}\label{lwl}\end{equation}
where $c, \tilde c$ are constants depending only on dimension. This is in the sense of expectation values; that is, for any semiclassical pseudodifferential operator $A$ on $\partial \Omega$ with principal symbol $a(y, \eta)$, $(y, \eta) \in T^* \partial \Omega$, we have \begin{equation}\begin{gathered}
\lim_{\lambda \to \infty}
\frac1{N_D(\lambda)}
\sum_{\lambda_j \le \lambda}
\lambda_j^{-2}\langle \psi_j, A_{\lambda^{-1}_j} \psi_j \rangle = c \int_{T^* (\partial \Omega)} (1 - |\eta|^2)^{1/2} 1_{\{ |\eta| \leq 1 \} } a(y, \eta) dy d\eta , \\
\lim_{\mu \to \infty} \frac1{N_N(\mu)}
\sum_{\mu_j \le \mu}
\ang{w_j, A_{\mu^{-1}_j} w_j} = c \int_{T^* (\partial \Omega)}  (1 - |\eta|^2)^{-1/2} 1_{\{ |\eta| \leq 1 \} } a(y, \eta) dy d\eta \end{gathered}\label{LWL}\end{equation}
where $N_D, N_N$ are the Dirichlet, resp. Neumann eigenvalue counting
functions. Here $y \in \partial \Omega$, $(y, \eta) \in T^* (\partial
\Omega)$ and $|\eta|^2$ is calculated with respect to the induced
metric on the boundary.  See \cite{hassell}.

Here we have adopted the semiclassical scaling, that is the
wavevectors 
at eigenvalue $\lambda_j^2$ are scaled by $h = h_j =
\lambda_j^{-1}$ or $\mu_j^{-1}$ so that they are rescaled to have length $1$ in
$\mathbb{R}^n$,
and therefore length $\leq 1$ when restricted to the boundary.

An intuitive explanation for the difference in the distribution of the
$\psi_j$ and the $w_j$ is as follows. If we use Fermi normal
coordinates $(y,r)$ near $\partial \Omega$, where $r$ is distance to
$\partial \Omega$, and if $(\eta, \rho)$ are the dual cotangent
coordinates, then the symbol of the semiclassical operator $h^2 \Delta
- 1$ is $\sigma(h^2 \Delta - 1) = \rho^2 + |\eta|^2 - 1$.  The
semiclassical normal derivative $ih d_n$ has symbol $\rho$, which when
restricted to the boundary and to the characteristic variety $\{
\sigma(h^2 \Delta - 1) = 0 \}$ is equal to $\sqrt{1 -
 |\eta|^2}$. Since the Dirichlet expectation value
$\lambda_j^{-2}\langle \psi_j, A_{\lambda^{-1}_j} \psi_j \rangle$
involves the application of two semiclassical normal derivatives (one
for each factor of $\psi_j$), compared to the Neumann expectation
value $\ang{w_j, A_{\mu^{-1}_j} w_j}$, it is not surprising that
the Dirichlet distribution in \eqref{LWL} is $1 - |\eta|^2$ times the Neumann
distribution. We can draw a moral from this.
\begin{equation}
\parbox{4.2in}{
\textbf{Moral: } semiclassically, in the Neumann case, the boundary trace that is analogous to  $d_n u$ in the Dirichlet case is not $v |_{\Omegab}$, but rather $(1 - h^2 \Deltab)_+^{1/2}  (v |_{\Omegab})$, where $\Deltab$ is the (positive) Laplacian on the boundary, $u, v$ are Helmholtz solutions at energy $h^{-2}$, and $(\dots)_+$
denotes the positive part.
}
\label{moral}\end{equation}
To see what goes wrong with using the naive measure $\tt$ of the `boundary condition error', let us attempt to  follow the same reasoning as in Section~\ref{sec:Dinc}. We can certainly  show that \be
(\min \tt[u])^{-2} \; = \;
\Bigl\|
\sum_j \frac{ w_j \ang{w_j, \cdot} }{(\Et - \mu_j^2)^2}
\Bigr\|_{L^2(\Omegab) \to L^2(\Omegab)}  .
\label{e:tIdinvsq}
\ee
The problem is that the $ w_j $ do not behave as uniformly as the Dirichlet traces $\psi_j$; we have a lower bound
\begin{equation}
\| w_j \|_{L^2(\partial \Omega)} \geq c,
\label{lower}\end{equation}
but the sharp upper bound is \begin{equation}
\| w_j \|_{L^2(\partial \Omega)} \leq C \mu_j^{1/3}.
\label{upper}\end{equation}
(This estimate follows from Tataru \cite{Tataru}.) The reason why, in Theorem~\ref{MPS-Dir}, we were able to get upper and lower bounds on $d(E, \spec_D)$ of the same order in $E$ was that the \emph{lower} bound on the operator norm of a \emph{single} term $\psi_j \ang{\psi_j, \cdot}$ was of the same order as the \emph{upper} bound on the sum $\sum_j  \psi_j \ang{\psi_j, \cdot}$ \emph{over a whole spectral cluster} $|\lambda - \lambda_j| \leq 1$. In the Neumann case, using $\tt$ will lead to a gap of at least $\mu^{1/3} = \Et^{1/6}$ between the upper and lower bounds on $d(\Et, \spec_N)$.
\vskip 8pt
Notice, however, that if we take our \textbf{Moral}, \eqref{moral}, seriously, then we could expect\break to find good upper and lower bounds on the quantity
$(1 - h_j^2 \Deltab)_+^{1/2} w_j$ instead. Indeed, this is the case, and we have the following exact analogues of  \eqref{upperlower}, and Theorem~\ref{main}:


\begin{theo}\label{Neu-est}
Let $\Omega \subset \RR^n$ be a smooth bounded domain, and let $w_j$ be the restriction to $\partial \Omega$ of the $j$th $L^2$-normalized Neumann eigenfunction $v_j$. Then there are constants $c, C$ such that \vskip 5pt
(i) $\|  (1 - h_j^2 \Deltab)_+^{1/2} w_j \|_{L^2(\partial \Omega)} \geq c$, \ $h_j = \mu_j^{-1}$;
\vskip 5pt
(ii) the operator norm of \begin{equation}
\sum_{\mu_j \in [\mu, \mu + 1]} (1 - h^2 \Deltab)_+^{1/2} w_j \ \big\langle(1 - h^2 \Deltab)_+^{1/2} w_j, \cdot \big\rangle, \quad h = \mu^{-1},
\label{Nopnorm}\end{equation}
is bounded by $C$.
\end{theo}

\begin{example} On the unit disc, Neumann eigenfunctions have the form $$
v(r, \theta) = ce^{in\theta} J_n(\mu_{n,l} r), \quad \mbox{ where } \;
J_n'(\mu_{n,l}) = 0, $$
and from \eqref{basic} we derive $$
2 \mu_{n,l}^2 = \int_{\partial \Omega} (\mu_{n,l}^2 - n^2) |v|^2 \implies \| (1 - \Deltab/\mu_j^2)_+^{1/2} w_j \| = \sqrt{2}. $$
Since zeroes of $J_n'$ are at least $\pi$ apart, we see that the operator norm  \eqref{Nopnorm} in the case of the unit disc is precisely $\sqrt{2}$. So we see in the case of the unit disc that Theorem~\ref{Neu-est} holds with $c = C = \sqrt{2}$.

Also note that when $l=1$,  $\mu_{n,1} \sim n + c n^{1/3}$, and then $\| w_j \| \sim \mu_j^{1/3}$. These are `whispering gallery modes', which saturate the bound \eqref{upper}. \end{example}

We now sketch some parts of the proof of Theorem~\ref{Neu-est}. Let us show how to obtain an upper bound for a single function $w_j$, that is, prove
$$
\|  (1 - h_j^2 \Deltab)_+^{1/2} w_j \|_{L^2(\partial \Omega)} \leq C.
$$
This already contains the crucial difficulties in the proof.

We return to \eqref{basic}, and deduce from it, using a vector field $V$ equal to $d_n$ at the boundary,  that
$$
\int_{\partial \Omega} v_j d_n^2 v_j = O(\mu_j^2).
$$
(This is not quite as straightforward as in the Dirichlet case, as one needs to show that  the left hand side is $O(\mu_j^2)$. This requires some integration-by-parts and relies on estimate \eqref{upper}.) It follows, using $(\Delta - \mu_j^2) v_j = 0$ at $\partial \Omega$, and that $\Delta = -d_n^2 + \Delta_{\partial \Omega}$ at the boundary, modulo first order operators, that $$
\int_{\partial \Omega} w_j ( (1 - h_j^2 \Deltab) w_j ) = O(1).
$$
That is,
$$
\| (1 - h_j^2 \Deltab)_+^{1/2} w_j \|_{L^2(\partial \Omega)}^2 - \| (h_j^2 \Deltab - 1)_+^{1/2} w_j \|_{L^2(\partial \Omega)}^2 = O(1). $$
So it remains to show that \begin{equation}
\| (h_j^2 \Deltab - 1)_+^{1/2} w_j \|_{L^2(\partial \Omega)}^2  = O(1) \quad \text{ (cf. \eqref{lwl}).}
\label{rts}\end{equation}

Intuitively, this quantity should be very small, since the operator $(h_j^2 \Deltab - 1)_+^{1/2}$ is microsupported in the elliptic region where we expect eigenfunctions should be negligible (the wavenumber on the boundary exceeds $\mu_j$
hence waves are evanescent in the normal direction).
However, there is a difficulty since the microsupport of $(h_j^2 \Deltab - 1)_+^{1/2}$ meets the boundary of the hyperbolic region $ \{ |\eta| \leq 1 \}$.

To prove \eqref{rts}, we break up $w_j$ spectrally into three pieces. Choose a smooth function $\phi$ of a real variable such that $\phi(t) = 0$ for $t \leq 5/4$ and $\phi(t) = 1$ for $t \geq 7/4$. Then we decompose 
$$ 
w_j = \big(1 - \phi \big)\big( \frac{h_j^2 \Deltab - 1}{h_j^{2/3}} \big) w_j + \phi\big( \frac{h_j^2 \Deltab - 1}{h_j^{2/3}} \big) \big(1 - \phi \big)\big(h_j^2 \Deltab - 1 \big) w_j + \phi \big(h_j^2 \Deltab - 1 \big) w_j,
$$
valid for $h_j$ sufficiently small, e.g. $h_j^{2/3} \leq 1/2$, and therefore 
\begin{equation}\begin{aligned}
(h_j^2 \Deltab - 1)_+^{1/2} w_j &= (h_j^2 \Deltab - 1)_+^{1/2} \big(1 - \phi \big)\big( \frac{h_j^2 \Deltab - 1}{h_j^{2/3}} \big) w_j \\ &+ (h_j^2 \Deltab - 1)_+^{1/2}  \phi\big( \frac{h_j^2 \Deltab - 1}{h_j^{2/3}} \big) \big(1 - \phi \big)\big(h_j^2 \Deltab - 1 \big) w_j \\ &+ (h_j^2 \Deltab - 1)_+^{1/2} \phi \big(h_j^2 \Deltab - 1 \big) w_j
\\
&:= I + II + III.
\end{aligned}\label{pieces}\end{equation}
Roughly speaking, here the first piece $I$ is supported in the frequency range $1\leq |\eta| \leq 1 + C h^{2/3}$, the second piece $II$ is supported in the frequency range $1 + c h^{2/3} \leq |\eta| \leq 2$ and the third piece $III$ is supported where $|\eta| \geq 3/2$. 
We now estimate each piece $I, II$ and $III$ separately. 

\emph{Estimating $I$}. We can estimate this piece purely using $L^2$ spectral theory. To do this we observe that  on the support of $(1 - \phi)(h^{-2/3}(h^2 \Deltab - 1))$, we have $(h^2 \Deltab - 1)_+^{1/2}  \leq 2h^{1/3}$. The required estimate now follows from this and \eqref{upper}. 

The other estimates are more intricate, and rely on expressing the boundary value of Neumann eigenfunctions in terms of themselves and the semiclassical double layer potential. More precisely, let $D_h$ denote the integral operator 
$$
D_h(x, y) = \frac1{2} \partial_{n_y} G_h(y,x), \quad x \neq y, 
$$
where $G_h(x,y)$ is the Helmholtz Green function $(\Delta - (h^{-1} + i0)^2)^{-1}(x,y)$ on $\RR^n$. It is well-known that 
$$
w_j = D_{h_j} w_j, \quad  h_j = \mu_j^{-1}.
$$
Iterating this we find that
\begin{equation}
w_j = D_{h_j}^{N} w_j, \quad N = 1, 2, \dots
\label{DN}\end{equation}
According to \cite{hassell},
$D_h^N$ is the sum of a semiclassical FIO microsupported in the hyperbolic region $\{ |\eta| < 1 \}$ in both variables; a pseudodifferential operator of order $-N$, in the sense that it maps $L^2(\Omegab)$ to $H^N(\Omegab)$ with norm $O(h^{N})$; and an operator microsupported close to $\{ |\eta| = 1 \}$ in both variables. 

\emph{Estimating $III$.} We use \eqref{DN} with $N=2$, and write $III$ as 
$$
(h_j^2 \Deltab - 1)_+^{1/2}\phi(h_j^2 \Deltab - 1) w_j = (h_j^2 \Deltab - 1)_+^{1/2}\phi(h_j^2 \Deltab - 1) D_{h_j}^2 w_j.
$$
When we compose $(h^2 \Deltab - 1)_+^{1/2}\phi(h^2 \Deltab - 1)$ with $D_h^2$, we obtain a pseudodifferential operator of order $-1$, that is, mapping $L^2(\Omegab)$ to $H^1(\Omegab)$ with norm $O(h)$, since $\phi(h^2 \Deltab - 1)$ is microsupported away from $\{ |\eta| \leq 1 \}$, and the estimate again follows from this and \eqref{upper}. 

\emph{Estimating $II$.} This is the most delicate estimate, since the spectral cutoff $\phi\big( h_j^2 \Deltab - 1 / h_j^{2/3} \big) $ is only localizing frequencies at a distance $h^{2/3}$ away from the hyperbolic set $\{ |\eta| \leq 1 \}$ where $D_h$ is an order zero operator. To deal with this we use the following lemma:

\begin{lem} The operator 
\begin{equation}
\phi\big( \frac{h^2 \Deltab - 1 }{ h^{2/3}} \big) (1 - \phi)\big( h^2 \Deltab - 1 \big) 
\label{sharpcutoff}\end{equation}
 can be represented as an oscillatory integral 
\begin{equation}
(2\pi h)^{-(n-1)} \int e^{i \Phi(y, y', \xi)/h } b \big( \frac{y+y'}{2}, \xi, h \big) \, d\xi,
\label{a}\end{equation}
where
\begin{equation}
\Phi(y, y', \xi) = \sum_{j,k=1}^{n-1} a_{jk} \big( \frac{y+y'}{2} \big) (y-y')_k \xi_j
\label{phase}\end{equation}
and $b$ satisfies estimates 
\begin{equation}
\Big| \partial^\alpha_y \partial^\beta_\xi b(y, \xi, h) \Big| \leq C_{\alpha \beta} h^{-2|\beta|/3}.
\label{symb}\end{equation}
\end{lem}

The proof will be given in a future article. 

\begin{remark} It does not seem possible to write the operator \eqref{sharpcutoff} in the usual pseudodifferential form, with phase function $(y-y') \cdot \xi$, because then the principal symbol would be 
$$
\phi \big( \frac{ g^{ij}(y) \xi_i \xi_j - 1}{h^{2/3}} \big).
$$
This function loses a factor $h^{-2/3}$ when differentiating in either $y$ or $\xi$ and such a symbol class does not lead to a sensible calculus (in the sense of having a composition formula, etc). By contrast, with a judicious choice of the $a_{jk}$ function in \eqref{phase}, one can arrange that the principal symbol of \eqref{a} is 
$$
\phi \big( \frac{ |\xi|^2 - 1}{h^{2/3}} \big),
$$
which satisfies the better estimates in \eqref{symb}, in that there is no loss of powers of $h$ when differentiating in $y$. 
\end{remark}

We can write down an oscillatory integral representation for the operator $D_h$ as an intersecting Lagrangian distribution. Using this and the  oscillatory integral representation for operator \eqref{sharpcutoff} given by the lemma, and using \eqref{DN} with $N=1$, we can write the operator  in $II$ as an oscillatory integral involving one factor of $D_h$. In this integral,  the phase is non-stationary on the support of the symbol. Using integration by parts in a standard way, we can show that this operator has an operator norm bound of $O(h^{1/3})$, which combined with \eqref{upper} gives the result.

\begin{remark} By using \eqref{DN} with large values of $N$, we can show that the $L^2$ norms of $II$ and $III$ are $O(h^\infty)$.
\end{remark}

\begin{figure} a)\raisebox{-1in}{\epsfig{file=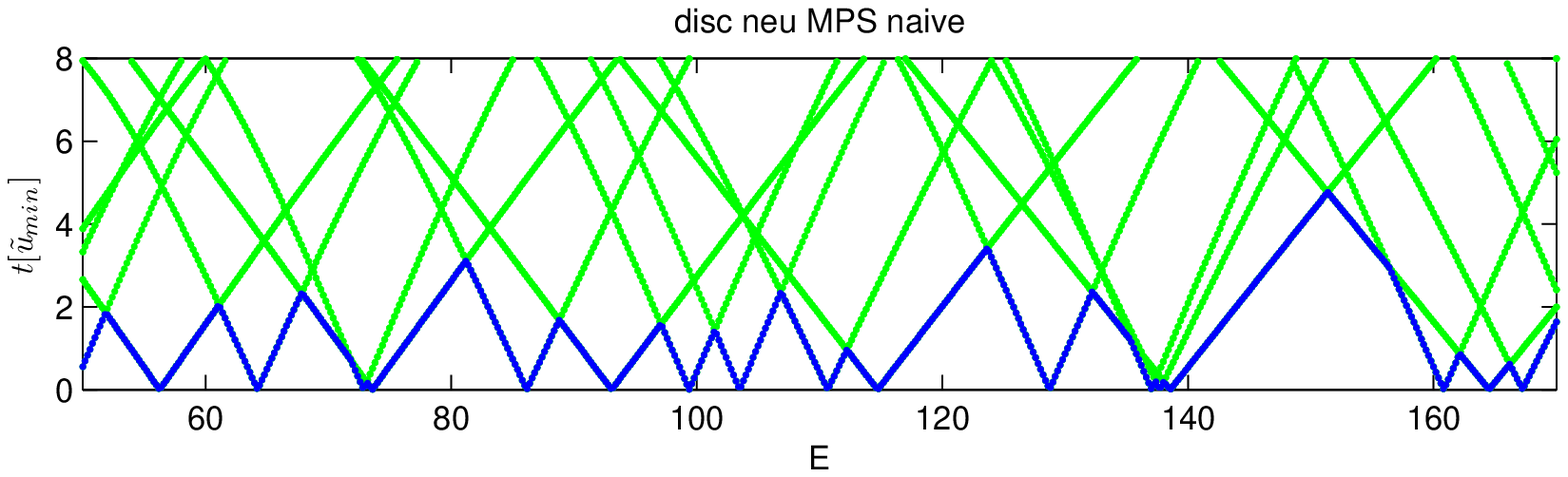,width=4in}}
b)\raisebox{-1in}{\epsfig{file=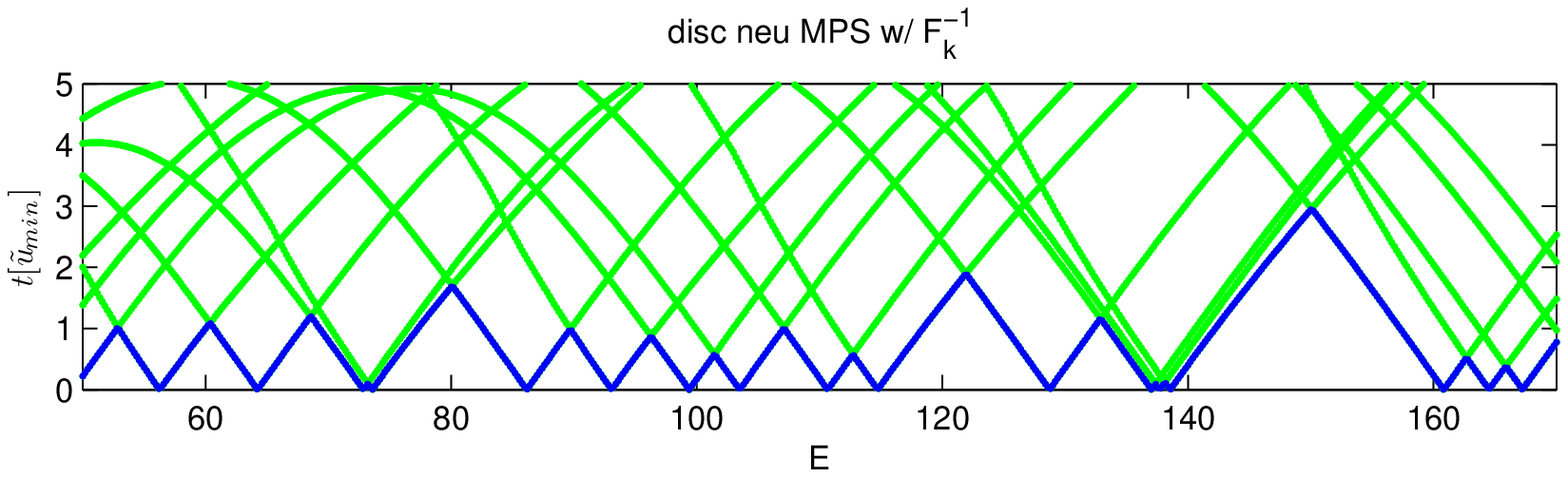,width=4in}}
\caption{
Graphs of minimum $\tttm[v]$ achievable at each $\tilde{E}$ for the
Neumann case,
for $\Omega$ the unit disc, and low energy.
a) $F=\mbox{Id}$, b) $F=F_\mu$ given by \eqref{e:Fmu}.
As in Fig.~\ref{f:tE}, $v$ is restricted to lie in a sufficiently
large numerical subspace.
The curves (shown lighter) lying above the lowest show higher
generalized eigenvalues of a matrix pair used to compute
$\tttm$ (the lowest gives $\tttm$ itself).
\label{f:tEneudisc}
}
\end{figure}

\begin{figure} a)\raisebox{-1.2in}{\epsfig{file=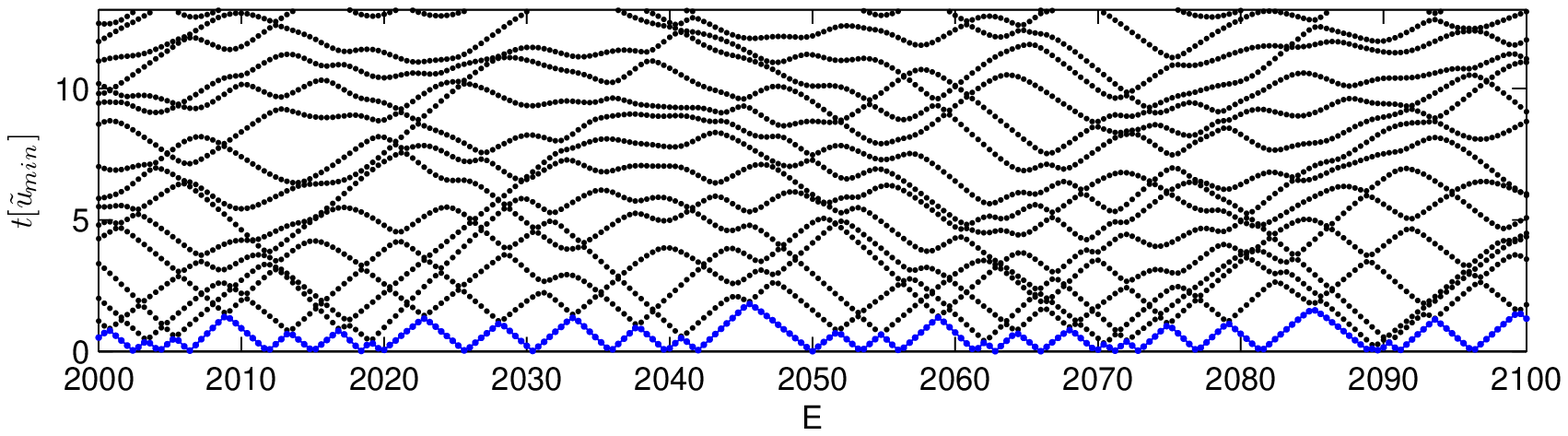,width=3.5in,height=1.3in}}
b)\raisebox{-1.2in}{\epsfig{file=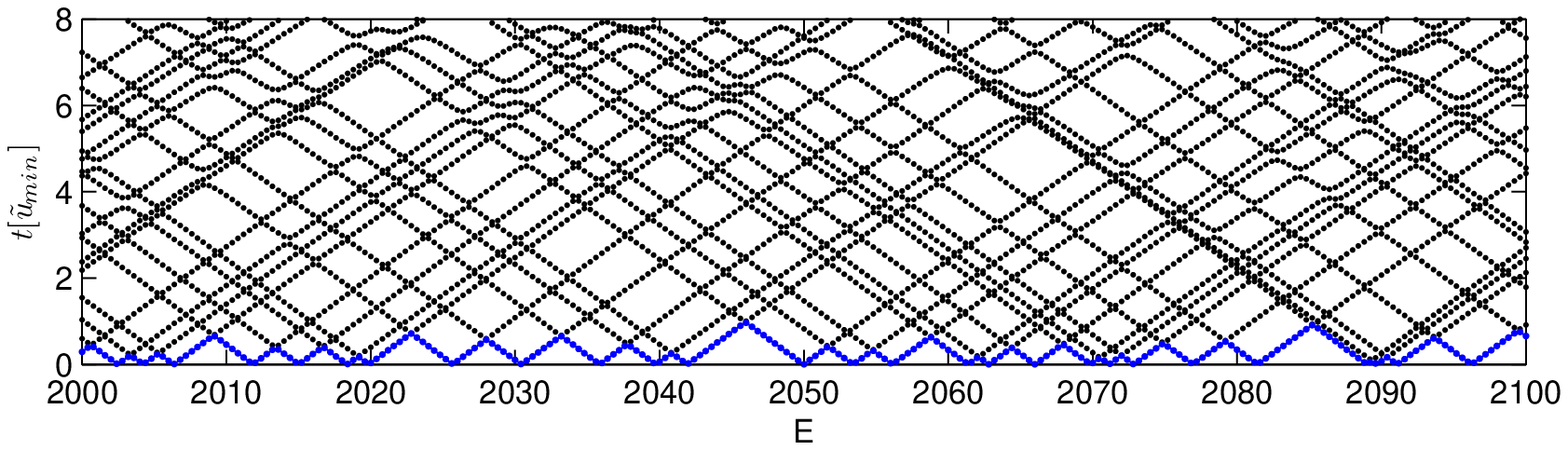,width=3.5in,height=1.3in}}
\caption{
Graphs of minimum $\tttm[v]$ achievable at each $\tilde{E}$ for the Neumann case,
for $\Omega$ the domain shown in Fig.~\ref{f:mode}, at higher energy.
Cases a) and b) are as in Fig.~\ref{f:tEneudisc}.
\label{f:tEneuhigh}
}
\end{figure}

\begin{figure} 
\epsfig{file=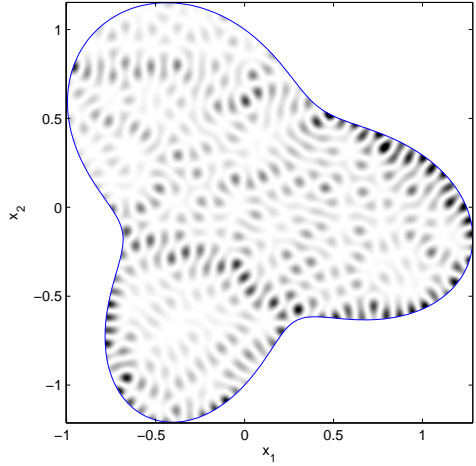,width=2.5in}
\caption{High-lying approximate Neumann eigenfunction of a smooth
domain $\Omega$
(density plot shows $|v(\mathbf{x})|^2$), computed by the new
MPS proposed in Section~\ref{s:neu}. $v$ minimizes
$\tttm[v]$ at  energy $\tilde{E}=2096.240170$ which is close to a local tension minimum.
The tension here was $\tttm < 10^{-6}$.
The eigenvalue is near the 500th.
Numerical computation used a basis of 400 fundamental solutions
lying outside $\Omega$ (see \texttt{MPSpack} manual \cite{polygonscatt}).
\label{f:mode}
}
\end{figure}

\section{Neumann inclusion bound and numerical demonstration}
\label{s:neu}

Using Theorem \ref{Neu-est} as a crucial tool, we propose the following
method of particular solutions (MPS)
for finding Neumann eigenpairs: at each energy $\tilde{E}=\mu^2$
we minimize the quantity
\begin{equation}
\tttm[v] = \frac{ \| F_\mu(\Deltab)(d_n v) \|_{L^2(\partial \Omega)} }{ \| v \|_{L^2(\Omega)} },
\label{ttt}\end{equation}
where (cf. \textbf{Moral} {\eqref{moral}})
the invertible boundary operator is $F_\mu(\Deltab)$ and
\be
F_\mu(\sigma) := \begin{cases}
\big( 1 - \mu^{-2} {\sigma} \big)^{-1/2}, \quad \sigma \leq \mu^2 - \mu^{4/3} \\
\mu^{1/3}, \phantom{aaaaaaaaaa} \quad \sigma \geq \mu^2 - \mu^{4/3}.
\end{cases}
\label{e:Fmu}
\ee
The effect of $F_\mu$ is roughly to boost the amplitudes of spatial frequencies
on the boundary which are close in magnitude to the overall wavenumber $\mu$;
however, it is regularized to limit this boost to a finite value
taking heed of the scaling \eqref{upper}.

This leads to the identity,
analogous to \eqref{e:tinvsq} and \eqref{e:tIdinvsq},
\begin{equation}
(\min_v \tttm[v])^{-2} = \Big\| \sum_j \frac{ F_\mu(\Deltab)^{-1} w_j \ang{F_\mu(\Deltab)^{-1} w_j, \cdot} }{(\mu^2 - \mu_j^2)^2} \Big\|,
\label{tttF}\end{equation}
where again, the min is to be taken over Helmholtz solutions $v$.
First we give a few words about why \eqref{tttF} holds. As before, $\displaystyle(\min_v \tttm[v])^{-1}$ is the operator norm of the composite function $g \mapsto f \mapsto v$, where $f = F_\mu(\Deltab)^{-1} g$ and $v$ is the Helmholtz solution with $d_n v = f$.
In a similar fashion to the derivation of \eqref{Zexp}, we have
$$\begin{gathered}
v = \sum_j \frac{\ang{f, w_j} v_j }{\mu^2 - \mu_j^2} = \sum_j \frac{\ang{F_\mu(\Deltab)^{-1} g, w_j} v_j }{\mu^2 - \mu_j^2} \\
= \sum_j \frac{\ang{ g, F_\mu(\Deltab)^{-1}w_j} v_j }{\mu^2 - \mu_j^2} .
\end{gathered}$$
Then a $T^* T$ argument, analogous to that in
the proof of Theorem~\ref{main}, gives \eqref{tttF}.

Finally, since $F_\mu(\Deltab)^{-1}$ is essentially $(1 - h^2 \Deltab)^{1/2}_+$, we can use Theorem~\ref{Neu-est} (together with \eqref{upper}) to prove the following tight Neumann inclusion bound (analogous to Theorem~\ref{MPS-Dir}):

\begin{theo}\label{MPS-Neu} There exist constants $c, C$ depending only on $\Omega$ such that the following holds.  Let $v$ be a nonzero solution of $(\Delta - \mu^2) v = 0$ in $C^\infty(\Omega)$. Let $\tttm[v]$ be as in \eqref{ttt}, and let $\vmin$ be the Helmholtz solution minimizing $\tttm[v]$. Then
$$
c   \tttm[\vmin] \leq d(\mu^2, \spec_N) \leq C \tttm[v].
$$
\end{theo}

We postpone the full proof to a future publication.
Comparing to the Dirichlet case, we note that there are no factors
of $\sqrt{\tilde{E}}$ in the bound (i.e. $\alpha=0$); this is to be expected dimensionally
since the Neumann tension
$\tttm[v]$ already contains an extra derivative compared to the Dirichlet
tension $t[u]$.

We end with some preliminary numerical demonstrations of our Neumann
MPS in $n=2$ dimensions.
The lowest curves shown in a) and b) of Fig.~\ref{f:tEneudisc} are the
Neumann analogs of Fig.~\ref{f:tE},
for the unit disc, comparing two choices of the operator $F$.
It is clear that the naive choice $F=\mbox{Id}$ leads to large variations
in slopes, whereas choosing $F = F_\mu$ given by \eqref{e:Fmu}
causes these slopes to become very similar.
As discussed, the disc allows Neumann modes whose value $L^2$ norms on the
boundary vary as widely as is possible.
Fig.~\ref{f:tEneuhigh} shows, at higher energies, the same
but for the smooth planar domain shown in Fig.~\ref{f:mode};
the difference in slopes of the lowest curve is less striking.
However, we have also plotted curves showing the higher generalized
eigenvalues relevant to the numerical implementation of the
MPS \cite{incl}. It is clear that our proposed $F$ operator
causes these higher curves to acquire not only very uniform slopes
but much less `interaction' between the curves, both of which
should lead to an improved numerical method.

Finally, in Fig.~\ref{f:mode}
we plot a Neumann eigenfunction computed with our proposed
MPS using $F$ as in \eqref{e:Fmu}. The tension $\tttm[v]$ was found to be
less than $10^{-6}$. The constant $C$ in Theorem~\ref{MPS-Neu}
is unknown, but the local slope of tension graph was
measured to be about 0.5, corresponding
to $C\approx 2$. Thus the inclusion bounds on the eigenvalue
are $[2096.240168, 2096.240172]$, i.e. about 9 digits of relative accuracy.
Computation took a few seconds on a laptop, using a basis set of size
400, and 450 quadrature points on $\partial \Omega$.
The $F$ operator was approximated to spectral accuracy
using trigonometric polynomials on the boundary.
We note that applying $F$ in higher dimensions $n\ge 3$ will prove
more of a challenge.

\bibliographystyle{amsplain} 
\bibliography{alex}

\end{document}